\documentclass[12pt]{article}

\usepackage{graphicx}

\usepackage{amssymb}
\usepackage{amsthm}
\usepackage{amsfonts}
\usepackage{amsmath}

\usepackage{alltt}

\usepackage{url}

\newtheorem{prop}{Proposition}
\newtheorem{lem}[prop]{Lemma}

\newtheorem{rem}{Remark}
\newtheorem{defn}{Definition}
\newtheorem{pf}{Proof}

\newcommand{\R}{\mathbb{R}}

\newcommand{\C}{\mathbb{C}}

\title{Fixed-Point Methods on Small-Signal Stability Analysis}

\author{Licio Hernanes Bezerra\footnote{Universidade Federal de Santa Catarina, Departamento de Matem\'atica,
Florian\'opolis, SC 88040-900, Brazil. email: licio.bezerra@ufsc.br}
}
\date{}

\begin{document}

\maketitle

\begin{abstract}
In this paper we introduce the Diagonal Dominant Pole Spectrum Eigensolver (DDPSE), which is a fixed-point method that computes
several eigenvalues of a matrix at a time.  DDPSE is a slight modification of the Dominant Pole Spectrum Eigensolver (DPSE), that has being used in
power system stability studies. We show that both methods have local quadratic convergence. Moreover, we present practical results obtained by both methods,
from which we can see that those methods really compute dominant poles of a transfer function
of the type $c^T(A-sI)^{-1}b$, where $b$ and $c$ are vectors, besides being also effective in finding low damped modes of a large scale power system.

\end{abstract}


\section{Introduction}


A power system can be described as a coupled system of differential and algebraic equations. The following system is obtained by linearizing the system at an operating point:
$$ \left( \begin{array}{c} \dot{x}  \\ 0 \end{array} \right) = \left( \begin{array}{cc} J_1 & J_2 \\ J_3 & J_4 \end{array} \right) \left( \begin{array}{c} x  \\ y \end{array} \right),$$
where $J_1,J_2,J_3,J_4$ are matrices, $x$ is the vector of dynamical variables and $y$ is the vector of algebraic ones. The matrix $J = \left(\begin{array}{cc} J_1 & J_2 \\ J_3 & J_4 \end{array}\right)$
denotes the Jacobian matrix of the system. Since $y = -J_4^{-1}J_3x$, we have
$$\dot{x} = \left( J_1 - J_2J_4^{-1}J_3 \right) x.$$
The matrix $A=J_1 - J_2J_4^{-1}J_3$ is called the state matrix of the system. For large scale power systems $J$ is very  sparse, while $A$ is not in general.
We can observe that, given a vector $x$, we can easily compute $z =(A-\lambda I)^{-1}x$
by solving the following algebraic system
$$ \left( \begin{array}{cc} J_1 & J_2 \\ J_3 & J_4 \end{array} \right) \left( \begin{array}{c} z  \\ w \end{array} \right) = \left( \begin{array}{c} x  \\ 0 \end{array} \right).$$
From that we note that standard methods of eigenvalue calculation can be used in order to compute eigenvalues of $A$, even without explicitly calculating the state matrix.

Knowledge of rightmost eigenvalues of $A$ is essential in the power system small-signal stability analysis. In the literature
there are several papers that use classical methods to compute rightmost eigenvalues of a state matrix from the above calculation \cite{an,du,is,nm1,rm}.
On the other hand, some authors prefer instead to deal with the generalized eigenvalue problem $J u = \lambda E u$, where
$$ E = \left( \begin{array}{cc} I & 0 \\ 0 & 0 \end{array} \right),$$
and $I$ is the identity matrix. However, this approach requires a non-obvious strategy to control instability caused by the spurious eigenvalue at infinity, for instance, if you use generalized M\"obius transforms \cite{lt,ll} .
The landscape of small-signal stability analysis has changed a little when methods based on transfer functions, like DPA \cite{nmleo} and DPSE \cite{nm}, have arisen in literature
\cite{nm,rog}. The Dominant Pole Algorithm (DPA), which computes a single eigenvalue at a time, is actually a Newton's method according to a simple calculation shown in \cite{lb}. The Dominant Pole Spectrum Eigensolver (DPSE), which can be seen as a generalization of DPA in a certain way, is a fixed-point method that can compute several eigenvalues at a time. On the one hand, each step requires the solution of $p$ linear systems if you want to calculate $p$ eigenvalues. On the other hand,
in power system stability studies, a suitable pre-ordering of the Jacobian matrix prevents large amounts of fill-in and thus
its sparse LU factorization is done with lower computational complexity.
Moreover, DPSE converges quadratically and a proof of its local quadratic convergence first appeared in \cite{lb}. Nevertheless, here we give an easier proof, which can be seen in \S 2.
We will also see that a slight modification of DPSE yields a new fixed-point method, the Diagonal Dominant-Pole Spectrum Eigensolver (DDPSE), that also has local quadratic convergence, as discussed in \S 3.
In the last section, we compare results obtained from the implementation of those two methods regarding the time of computation. From these tests we verify that both methods really compute dominant poles of a transfer function
of the type $c^T(A-sI)^{-1}b$, where $b$ and $c$ are vectors, besides being also effective in finding low damped modes of the system.


\begin{figure}[htbp]
\centering
		\includegraphics[width=1.0\textwidth]{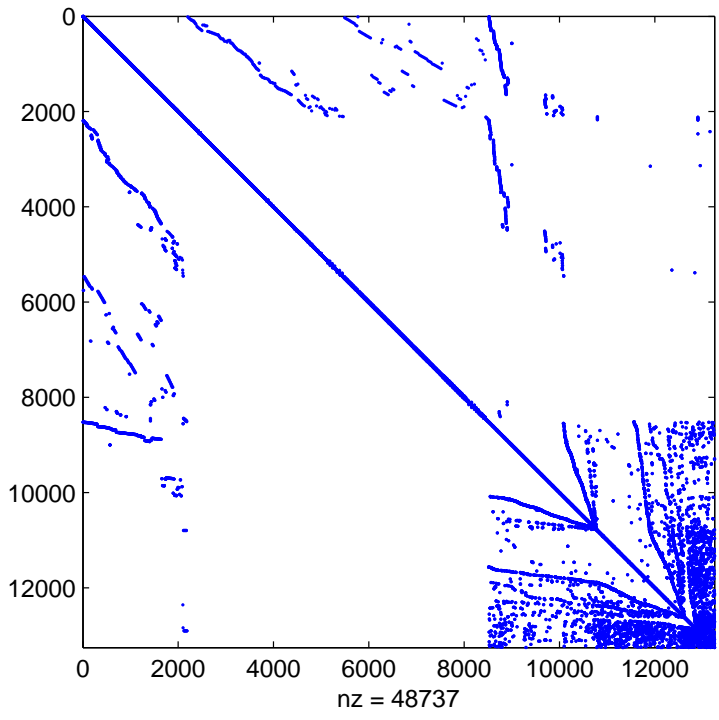}
\caption{Sparse pattern of $M$ (nz is the number of its nonzero entries)}
\end{figure}

\section{DPSE}

The motivation for the Dominant Pole Spectrum Eigensolver (DPSE) came from SISO dynamical systems $(E,J,B,C,D)$
of the form
$$\left\{ \begin{array}{rcl}
E \dot{x} (t) & = & J x(t) +  Bu(t) \\
y(t) & = & C^T  x(t) + D u(t)
\end{array}
\right.
$$
where $J,E \in \R^{N\times N}$, $E=diag([1,...,1,0,...,0])$;
$x(t) \in \R^{N\times 1}$
is composed by dynamical and algebraic variables, $x_d(t)$ and $x_a(t)$, which are respectively associated with the unit and the null diagonal entries of $E$;  $B, C \in \R^{N\times 1}$, where $B^T = \left( B_d^T \, B_a^T \right)$ and $C^T=\left( C_d^T \, C_a^T \right)$; $u(t)\in \R$ is the input, $y(t)\in \R$
is the output, and $D\in \R$. The corresponding
transfer function $h:\C \to \C$ is defined as
$$ h(s) = C^T \left( sE - J\right)^{-1} B + D.$$
Suppose that there are $n$ dynamical variables in the system. If $J_1=J(1:n,1:n)$, $J_2 = J(1:n,n+1:N)$, $J_3=J(n+1:N,1:n)$
and $J_4=J(n+1:N,n+1:N)$, then
\begin{equation}\label{equal}
h(s) = c^T \left( sI- \left(J_1-J_2J_4^{-1}J_3\right)\right)^{-1}  b + d,
\end{equation}
where $b=B_d - J_2J_4^{-1} B_a$, $c=C_d - J_3^TJ_4^{-T} C_a$ and $d = D - C_a^TJ_4^{-1}B_a$.
The matrix $A=J_1-J_2J_4^{-1}J_3$ is called the state matrix of the system.
Note that, for any $\mu \notin \lambda (A)$ and for any $b\in \C^n$,
$$\left( \begin{array}{c} (A-\mu I)^{-1}b \\ 0
\end{array}\right) =E\, (J-\mu E)^{-1} \left( \begin{array}{c} b\\0
\end{array}\right).$$

Suppose that $A\in \R^{n\times n}$ is diagonalizable, that is,
$A=PDP^{-1}$, where $P$ is an invertible matrix and $D$, a diagonal matrix.
So, the spectrum of $A$, denoted by $\lambda (A)$, is the set of the diagonal entries of $D$.
From now on we suppose that every eigenvalue of $A$ is simple.

Let $b,c\in \R^n$ such that $c^T (A-s I)^{-1} b\ne 0$ for all $s\in \C - \lambda(A)$,
and $c^TPe_k e_k^TP^{-1}b \ne 0$ for $k=1:n$.
If $d=0$, from Equation~\ref{equal},
$$ h(s) = \sum_{k=1}^{n} \frac{R_k}{d_{kk}-s},$$
where $R_k = c^TPe_k e_k^TP^{-1}b$, $k=1:n$.
\begin{defn}
A pole $d_{kk}$ is called a dominant pole if it corresponds to a relatively large
$m_k = \dfrac{|R_k|}{|Re(d_{kk})|}$. $m_k$ is the measure for dominance of the pole $d_{kk}$.
\end{defn}

Now, since
$$\frac{ (A-s I)^{-1} b}{c^T (A-s I)^{-1} b}= \frac{Adj \, (A-s I) b}{c^T Adj\, (A-s I) b} \mbox{ and }
\frac{ (A^T-s I)^{-1}c }{c^T (A-s I)^{-1} b}= \frac{Adj\, (A^T-s I)c}{c^T Adj\, (A-s I)b},$$
we conclude that the functions $f: \C \to \C^n$  and $g: \C \to \C^n$ defined respectively by
\begin{equation}\label{f(s)}
f(s) = \left\{ \begin{array}{cl} \frac{ (A-s I)^{-1} b}{c^T (A-s I)^{-1} b} & \mbox{ for } s\in \C - \lambda(A);\\
\frac{Pe_j}{c^TPe_j} & \mbox{ for } s=d_{jj}, j=1:n. \end{array} \right.
\end{equation}
and
\begin{equation}\label{g(s)}
g(s) = \left\{ \begin{array}{cl} \frac{ (A^T-s I)^{-1} c}{c^T (A-s I)^{-1} b} & \mbox{ for } s\in \C - \lambda(A);\\
\frac{P^{-T}e_j}{b^TP^{-T}e_j} & \mbox{ for } s=d_{jj}, j=1:n. \end{array} \right.
\end{equation}
are entire functions (bear in mind that any entry of the Classical Adjoint of $(A-sI)$ is a sum of products of its elements).

Let $S=\left( s_1\, ...\, s_p\right)^T \in \C^p$, $p\le n$, and suppose that $X(S),Y(S) \in C^{n\times p}$ are defined by
$X(S)e_k = f(s_k) \mbox{ and } Y(S)e_k= g(s_k), \mbox{ for } k=1:p,$ where $e_1,...,e_n$ are the canonical vectors.

\begin{lem} Let $S_0=\left( d_{k_1,k_1}\, ...\, d_{k_p,k_p}\right)$, where $ d_{k_1,k_1}\, ...\, d_{k_p,k_p}$ is a $p$-uple of distinct eigenvalues, $1\le k_1<...<k_p \le n$.
Then, there is an open neighborhood ${\cal O}$ of $S_0$ so that
$Y(S)^TX(S)$ is invertible for any $S$ belonging to ${\cal O}$.
\end{lem}
\begin{pf}
The lemma follows because
\begin{equation}\label{ytx-1}
Y(S_0)^TX(S_0)= {\rm diag} \left(\left[ \frac{1}{c^TPe_{k_1}e_{k_1}^TP^{-1}b},...,
\frac{1}{c^TPe_{k_p}e_{k_p}^TP^{-1}b}\right]\right),
\end{equation}
that is, $Y(S_0)^TX(S_0)$ is invertible.
\qed\end{pf}

Let $F: {\cal O} \to \C^{n\times n}$ defined by $F(S) = \left( Y(S)^TX(S)\right)^{-1} \left( Y(S)^TAX(S)\right)$. We see that $F$ is analytic. Since
$$F\left( \left( d_{k_1,k_1}\, ...\, d_{k_p,k_p}\right)^T \right) = {\rm diag} \left(\left[ d_{k_1,k_1},...,d_{k_p,k_p}\right]\right),$$ every eigenvalue of $F$ is simple for $S$ belonging to
an open subset $ {\cal Z}$ of ${\cal O}$.
Let $G: {\cal Z}\to \C^p$ be the function defined by $G(S) = \left( \lambda_1(F(S)) \, ... \, \lambda_p(F(S))\right)$,
where $\lambda_1(F(S)) < ... < \lambda_p(F(S))$ are the eigenvalues of $F(S)$ (for some order on the complex numbers).
Observe that $\left( d_{k_1,k_1}\, ...\, d_{k_p,k_p}\right)^T$ is a fixed point of $G$.
On the other hand, if $s_i$ is not an eigenvalue for $i=1:p$, we conclude that $F(S)$ is equal to
$${\rm diag} (S)
 + \left( Y(S)^TX(S)\right)^{-1}e e^T {\rm diag}\left(\left[ \frac{1}{c^T(A-s_1I)^{-1}b},...,
 \frac{1}{c^T(A-s_pI)^{-1}b}\right]\right),$$
for $Y(S)b=e$, where $e=ones(p,1)$.
In order to calculate the derivative of $G$, we first see that
$$ \frac{\partial F}{\partial s_i} \left( d_{k_1,k_1},...,d_{k_p,k_p}\right) = e_ie_i^T - \frac{1}{c^TPe_{k_i}e_{k_i}^TP^{-1}b}ve_i^T,$$
where $v^T = \left( c^TPe_{k_1}e_{k_1}^TP^{-1}b\, ... \, c^TPe_{k_p}e_{k_p}^TP^{-1}b\right)$.
Therefore,
$$ \frac{\partial \lambda_k}{\partial s_i} \left( d_{k_1,k_1},...,d_{k_p,k_p}\right)=
\langle \, \frac{\partial \lambda_k}{\partial a_{rs}} \left( F\left( d_{k_1,k_1},...,d_{k_p,k_p}\right)\right)\, , \,\frac{\partial F}{\partial s_i} \left( d_{k_1,k_1},...,d_{k_p,k_p}\right) \,\rangle.$$
Notice that
$$\frac{\partial \lambda_k}{\partial a_{rs}} \left( F\left( d_{k_1,k_1},...,d_{k_p,k_p}\right)\right) = e_k e_k^T.$$
Hence,
$$ \frac{\partial \lambda_k}{\partial s_i} \left( d_{k_1,k_1},...,d_{k_p,k_p}\right)= e_i^Te_k -\frac{1}{c^TPe_{k_i}e_{k_i}^TP^{-1}b}e_i^T e_ke_k^T v=0.$$

\begin{defn}[DPSE] The fixed-point iteration applied to the function
$$G(S) = \left( \lambda_1(F(S)) \, ... \, \lambda_p(F(S))\right),$$ where $F(S) = \left( Y(S)^TX(S)\right)^{-1} \left( Y(S)^TAX(S)\right),$
defines the Dominant Pole Spectrum Eigensolver.
\end{defn}

With this definition, we have just proved the following proposition:

\begin{prop}[DPSE converges at least quadratically]
Let $\lambda_1$, ..., $\lambda_p$ be $p$ distinct eigenvalues of $A$. Then, there is a neighborhood $V$ of $\left( \lambda_1\, ...\, \lambda_p\right)^T$ such that  DPSE converges at least quadratically to $\left( \lambda_1 \, ... \, \lambda_p\right)^T$ for any $S_0\in V$.
\end{prop}

\begin{rem}
Suppose that DPSE has just calculated $S^{(r)}=\left(\, s_1^{(r)}\, ... \, s_p^{(r)}\right)^T$ from
$S^{(r-1)}$, and
$X(S^{(r)})$ has not been computed  yet.
For any $k=1:p$, we have that
$$AX(S^{(r-1)})e_k= $$
$$ =(A-s_k^{(r-1)}I+s_k^{(r-1)}I)\,
\frac{(A-s_k^{(r-1)})^{-1}b}{c^T(A-s_k^{(r-1)})^{-1}b} =$$
$$=\frac{b}{c^T(A-s_k^{(r-1)})^{-1}b} + s_k^{(r-1)}X(S^{(r-1)})e_k.$$
Therefore, 
the relative error, $|| (A - s_k^{(r)}) X(S^{(r-1)})e_k|| / || X(S^{(r-1)})e_k||$,
becomes as follows:
\begin{equation}\label{stop}
\frac{|| \frac{b}{c^T(A-s_k^{(r-1)})^{-1}b} + (s_k^{(r-1)}-s_k^{(r)})\, X(S^{(r-1)})e_k||}{||X(S^{(r-1)})e_k||}.
\end{equation}
Notice that, if $s_k^{(r-1)}$ tends to an eigenvalue, then $c^T(A-s_k^{(r-1)})^{-1}b$ tends to zero, and so,
$X(S^{(r-1)})e_k$ is an approximation of an associated eigenvector.
\end{rem}

\begin{rem}
Suppose that, for some $j\in\{ 1,...n\}$, $d_{jj}$ is a converged eigenvalue at step $r$. Then, the corresponding right and left vectors, $x_r$ and $y_r$, are such that
$$y_r^Tx_r \approx \frac{e_j^TP^{-1}}{e_j^TP^{-1}b} \frac{Pe_j}{c^TPe_j} = \frac{1}{c^TPe_je_j^TP^{-1}b} = \frac{1}{R_j}.$$
\end{rem}

\section{DDPSE}

\begin{figure}[htbp]
	\includegraphics[width=\textwidth]{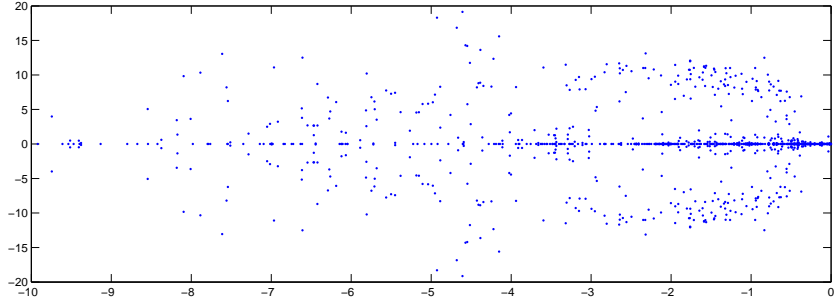}
\caption{Partial spectrum of the state matrix}
\end{figure}

If the fixed-point method is now applied to the diagonal of the matrix $F(S)$, then we have a variation of the DPSE, which will be called here as the Diagonal Dominant Pole Spectrum Eigensolver (DDPSE):
$$H\left(\left( s_1\, ...\, s_p\right)^T\right) \stackrel{\rm def}{=} {\rm diag}\left( F(s_1\, ... \, s_p)\right) =(s_1\, ...\, s_p)^T + $$
$$+ {\rm diag}\left( \left[ \frac{1}{c^T(A-s_1I)^{-1}b},...,\frac{1}{c^T(A-s_pI)^{-1}b}\right]\right) \left( Y(S)^TX(S)\right)^{-1} e$$
where $e=ones(p,1)$ and $s_i\notin \lambda(A)$ for all $i=1:p$.
Note that for $i=1:p$
$$\frac{1}{c^T(A-s_iI)^{-1}b}= \frac{\det (A-s_i I)} {c^T Adj(A-s_i I) b},$$
which is zero if $s_i$ is an eigenvalue of $A$.
So, we define $$H\left( (\lambda_1\, ...\, \lambda_p)^T\right) = (\lambda_1\, ...\, \lambda_p)^T$$
if $(\lambda_1\, ...\, \lambda_p)^T$ is a $p$-uple of distinct eigenvalues. Therefore, any $p$-uple of distinct eigenvalues is a fixed-point of $H$.

\begin{prop}[DDPSE converges at least quadratically] Let $\lambda_1$, $\lambda_2$, ..., $\lambda_p$
 be $p$ distinct eigenvalues of $A$. Then, there is a neighborhood $V$ of $(\lambda_1\, ...\, \lambda_p)^T$ such that, given any $S_0\in V$,  DDPSE converges at least quadratically to $(\lambda_1\, ..., \, \lambda_p)^T$.
\end{prop}
\begin{pf}
Let $\lambda_1=d_{k_1,k_1}$, ..., $\lambda_p= d_{k_p,k_p}$. Then,
$$ \frac{\partial H}{\partial s_i} \left( (d_{k_1,k_1}\, ...\, d_{k_p,k_p})\right) = e_i - \frac{1}{c^TPe_{k_i}e_{k_i}^TP^{-1}b}e_ie_i^T \left( Y(S)^TX(S)\right)^{-1} e.$$
From \ref{ytx-1}, $\left( Y(S)^TX(S)\right)^{-1} = {\rm diag} \left( \left[ c^TPe_{k_1}e_{k_1}^TP^{-1}b,...,c^TPe_{k_p}e_{k_p}^TP^{-1}b\right]\right)$. So,
$$ \frac{\partial H}{\partial s_i} \left( (d_{k_1,k_1}\, ...\, d_{k_p,k_p})\right) =0.\qed$$
\end{pf}

\begin{table}[htbp]

\begin{verbatim}
    k    ITER       CPU         k    ITER      CPU

   16     21      21.9805      19      3      2.8080
    4     22      23.0569      20      3      2.8080
   17     22      23.0881      16      4      3.7128
   20     23      23.9930      18      4      3.7128
    3     24      24.7418       1      5      4.3992
    7     24      24.7418      17      5      4.4304
    8     24      24.7418      14      6      5.0700
   12     24      24.7418      15      6      5.0700
   19     24      24.7730      13      7      5.6940
    5     25      25.3658       8      9      6.8016
    6     25      25.3658       9      9      6.8016
    9     25      25.3658      11      9      6.8016
   11     25      25.3658       2     10      7.1448
   14     25      25.3658      10     10      7.1760
   18     25      25.3970      12     10      7.1760
   15     26      25.8026       6     11      7.4880
    2     28      26.4422       7     11      7.5192
    1     33      27.5654       5     12      7.7532
   13     34      27.7838       4     15      8.2993
   10     35      27.8930       3     17      8.4709

\end{verbatim}

\caption{DDPSE versus DPSE}

\end{table}

\section{Numerical results}

\begin{table}[htbp]
{\small
\begin{alltt}
         Re      Imag        \(m\sb{k}\)          Re       Imag         \(m\sb{k}\)

      -0.6120   0.3587      12.40      -0.0335   -1.0787     760.11
      -2.9957  -9.3891       0.57      -5.7475    6.7761       1.43
      -4.5931  -0.2765       0.71      -1.2786    7.2546       1.97
      -7.5416  -6.2291       1.07      -4.3601    0.9544       0.48
      -1.2891  -8.5414 	     2.41      -4.1103    0.4801       1.04
      -2.9445  -4.8214 	     6.85      -1.2936    1.4028       5.59
      -4.0233   4.2124 	     2.61      -3.1928    9.2818       1.38
      -2.9445   4.8214 	     6.85      -1.8415    6.9859       5.11
      -9.7433   3.9765       0.14      -0.6120    0.3587      12.40
      -2.2927   0.0000 	     2.92      -2.9445    4.8214       6.85
      -5.8148   4.8704 	     1.36      -0.5776    6.2565       0.29
      -3.1928   9.2818       1.38      -0.5550    4.1191       0.35
      -0.0335   1.0787     760.11      -0.6786    7.1071       0.36
      -0.5567  -3.6097 	    14.87      -1.2891    8.5414       2.41
      -0.0335  -1.0787 	   760.15      -1.4790    8.2551       3.68
      -0.9401   8.1931 	     0.37      -0.5208    2.8814       0.78
      -1.2786   7.2546 	     1.96      -0.0335    1.0787     760.11
      -5.7475   6.7761 	     1.43      -0.5567    3.6097      14.87
      -5.5632   7.7510 	     1.22      -0.7584    4.9367       5.11
      -1.4790   8.2551 	     3.68      -0.4548    4.7054       5.78

\end{alltt}
}
\caption{Eigenvalues and their measure of relative dominance calculated by DDPSE and DPSE}
\end{table}

Our test matrix $J$ is sparse (density about 0.028\%), of order $N=13251$.
The pencil $Jv = \lambda E v$, where $E=diag(1,...,1,0,...,0)$, corresponds to
the problem $Ax = \lambda x$, where $A=J_1 - J_2J_4^{-1}J_3$
is of order $n = 1664$.
This is the Jacobian matrix that corresponds to a planning model of the Brazilian Interconnected Power System
and that had already used for tests in \cite{rog}.
In the tests with DPSE and DDPSE, we have used data which can be obtained from a specific transfer function. From this, $D=0$,  and the input vector $B= \left( B_d^T \, B_a^T\right)$ and the output vector $C=\left( C_d^T \, C_a^T\right)$ are as follows: $B(524)= B(1442) = 1$, $B(1884)=B(1918)=-1$, and
the others entries of $B$ are null; $C(11558)= 26.5721$, $C(11559)=-13.1127$, $C(12502)=-29.2954$, $C(12503)=3.7609$ and
these are all the non-zero elements of $C$. By using MATLAB, $C_a^T J_4^{-1} B_a=0$. Note that, from Equation~\ref{equal}, for any $s\notin \lambda(A)$, $C^T (J-sE)^{-1}B = c^T(A-sI)^{-1}b$,
where $b=B_d - J_2J_4^{-1} B_a$, $c=C_d - J_3^TJ_4^{-T} C_a$.

Let $s^{(k)} = (s_1^{(k)} \, ... \, s_p^{(k)})^T \in \C^p$. Let $X=X^{(k)}$ and $Y=Y^{(k)}$ be two matrices $N\times p$, such that for $j=1:p$ $X(:,j) = (J - s_j^{(k)}E)^{-1} B / C^T(J-s_j^{(k)}E)^{-1}B$ and
$Y(:,j) = (J^T - s_j^{(k)}E)^{-1} C / C^T(J-s_j^{(k)}E)^{-1}B$.
Let $V=V^{(k)}$ and $W=W^{(k)}$ be two matrices $n\times p$ such that, for $j=1:p$,
$$V(:,j) = (A - s_j^{(k)}I)^{-1} b / c^T(A-s_j^{(k)}I)^{-1}b$$ and
$$W(:,j) = (A^T - s_j^{(k)}I)^{-1} c / c^T(A-s_j^{(k)}I)^{-1}b.$$ So, we have
$$ W^TV = Y^T E X,$$ and for $j=1:p$
$$ W^TAVe_j = W^T (A - s_j^{(k)}I + s_j^{(k)}I) (A - s_j^{(k)}I)^{-1} b / c^T(A-s_j^{(k)}I)^{-1}b = $$
$$ = W^T b / c^T(A-s_j^{(k)}I)^{-1}b +  s_j^{(k)} W^T Ve_j = e / c^T(A-s_j^{(k)}I)^{-1}b + s_j^{(k)} W^T Ve_j,$$
where $e = ones(p,1)$. Therefore,
$$ F = (W^TV)^{-1} W^TAV = (W^TV)^{-1} ev^T+ S,$$
where $S=S^{(k)} = diag( s^{(k)} )$, and
$$v^T = \left( v^{(k)}\right)^T= \left( 1/ c^T(A-s_1^{(k)}I)^{-1}b \, ... \, 1/ c^T(A-s_p^{(k)}I)^{-1}b \right).$$
Hence, in order to use DPSE with $A$, $b$ and $c$,
we can carry out all the computation with $J$, $B$ and $C$ without explicitly computing $A$, $b$ and $c$.

We have specified a relative error tolerance of $10^{-5}$ to both right and left vectors. Here we have used
$$ \frac{||(J - s_j^{(k)}E) X^{(k-1)} e_j||}{||X^{(k-1)} e_j||} \mbox{ and }
\frac{||(J^T - s_j^{(k)}E) Y^{(k-1)} e_j||}{||Y^{(k-1)} e_j||}
$$
as the relative errors instead of using the equivalent formulae obtained from Equation~\ref{stop}, for we do not want to compute $b$ explicitly.
Convergence is achieved when both the convergence criteria are satisfied.
Suppose we have just obtained the first converged value $\lambda_1=s_1^{(k)}$. Then, we save their corresponding right and left generalized eigenvectors, $X^{(k-1)}e_1$ and $Y^{(k-1)}e_1$, to the respective columns of $X^{(r)}$ and $Y^{(r)}$ for $r\ge k$. Notice that the respective columns of $V^{(r)}$ and $W^{(r)}$ are formed by the dynamical variables of $X^{(k-1)}e_1$ and $Y^{(k-1)}e_1$, respectively. That eigenvalue can be deflated from the problem by this procedure.
To see that,
suppose $Ve_2 \approx a_1 v_1 + a_2 v_2$, where $v_1\approx V^{(k-1)}e_1$ and $v_2$ is an eigenvector corresponding to $\lambda_2$.
Hence,
$$A Ve_2 \approx a_1 \lambda_1v_1 + a_2 \lambda_2v_2\approx a_1 (\lambda_1-\lambda_2)Ve_1 + \lambda_2Ve_2,$$
and so,
$$\left( W^TV\right)^{-1} W^TA Ve_2 \approx a_1 (\lambda_1-\lambda_2)e_1 + \lambda_2e_2.$$
Thus,
$$F = \left( \begin{array}{ccccc} \lambda_1 & \times & \times &\cdots & \times \\
                                           0     & \lambda_2 & \times & \cdots & \times \\
                                           0     &   0       & \times & \cdots & \times \\
                                       \vdots    &  \vdots   & \vdots  & \cdots &\vdots \\
                                          0      &   0       & \times & \cdots  & \times
\end{array}\right).$$


In Table 1, for each eigenvalue calculated by DDPSE, or by DPSE, we list the CPU time (in seconds) and the number of iterations required for convergence,
together with the corresponding measure $m_k$ for dominance of a pole that was computed by the respective method.
For the tests we have started DPSE and DDPSE by choosing the same 20 initial shifts given in the left complex semi-plane: $\mu_k= k * (-1/20 + i/2)$, $k=1:20$.

The tests have been performed in the MATLAB R2011b 64 bits at a HP Compaq 6000 Pro, with processor Intel Core 2 Duo E8400 3.00 GHz.

\subsection{Conclusions}

In the tests, DPSE has shown that it is more stable than DDPSE. For instance, both algorithms started with the same 20 complex numbers in the upper-half plane. However, while DPSE converged to
19 eigenvalues still located in the upper half-plane, DDPSE converged to only 13 eigenvalues with positive imaginary part.
On the other hand, both converged to $-0.0335  \pm 1.0787\, i$, which are the most dominant poles of the system.
The DDPSE algorithm typically converges more slowly than the DPSE in total computer time, and we see that in Table 1.
Note that the eigenvalues of $A$ are clustered around zero, according to Figure 2, and even so both algorithms converged to dominant poles.



\begin{thebibliography}{99}

\bibitem{an} {\sc G. Angelidis and A. Semlyen}, {\it Improved methodologies for the calculation
of critical eigenvalues in small signal stability analysis}, IEEE
Trans. Power Syst., vol. 11, no. 3, pp. 1209--1217, 1996.

\bibitem{lb}
{\sc L. H. Bezerra}, {\it An eigenvalue method for calculating dominant poles of a transfer function}, Appl. Math. Lett. 21(3), pp. 244--247, 2008.

\bibitem{lt}
{\sc L. H. Bezerra and C. Tomei}, {\it Spectral transformation algorithms for
computing unstable modes}, Comp. Appl. Math, vol. 18, no. 1, pp.
1--15, 1999.


\bibitem{du} {\sc Z. Du, W. Liu, and W. Fang}, {\it Calculation of rightmost eigenvalues in
power systems using the Jacobi-Davidson method}, IEEE Trans. Power
Syst., vol. 21, no. 1, pp. 234--239, 2006.


\bibitem{is}
{\sc F. Ishikawa, H. Sasaki, J. Kubokawa, and H. Terasako}, {\it An efficient algorithm
for small signal stability analysis using the Arnoldi and S-matrix
methods}, Elect. Eng. Jpn., vol. 121, no. 4, pp. 38--47, 1997.

\bibitem{ll}
{\sc  L. T. G. Lima, L. H. Bezerra, C. Tomei, and N. Martins}, {\it New methods
for fast small-signal stability assessment of large-scale power systems},
IEEE Trans. Power Syst., vol. 10, no. 4, pp. 1979--1985, Nov. 1995.

\bibitem{nm1}
{\sc N. Martins}, {\it Efficient eigenvalue and frequency response methods applied to power system small-signal stability studies},
IEEE Trans. Power Syst. vol. 1, no. 1, pp. 217--226, 1986.

\bibitem{nmleo} {\sc N. Martins, L. T. G. Lima, H. J. P. Pinto}, {\it Computing Dominant Poles of Power System Transfer Functions},
IEEE Trans. Power Syst. vol. 11, no. 1, pp. 162--170, 1996.

\bibitem{nm}
{\sc N. Martins}, {\it The Dominant Pole Spectrum Eigensolver},
IEEE Trans. Power Syst. vol. 12, no. 1, pp. 245--254, 1997.



\bibitem{rm}
{\sc J. Rommes, N. Martins and F. Freitas}, {\it Computing Rightmost Eigenvalues for  Small-Signal Stability Assessment of Large-Scale Power Systems}, IEEE Trans. Power Syst. vol. 25, no. 2, pp. 929--938, 2010.


\bibitem{rog}
{\sc J. Rommes and G. L. G. Sleijpen}, {\it Convergence of the dominant pole algorithm and Rayleigh quotient iteration},
SIAM J. Matrix Anal. Appl., vol. 30, no. 1, pp. 346--363, 2008.


\end{thebibliography}
\end{document}